\documentclass{amsart}
\usepackage{amsmath}
\usepackage{amssymb}
\usepackage{enumerate}
\usepackage{enumitem}
\usepackage{mathrsfs}
\usepackage{cases}
\usepackage{color}
\usepackage{tikz}
\usepackage{appendix}
\usepackage{makecell}
\usepackage{multirow}
\usepackage{algorithmic,algorithm}
\usepackage{graphicx}
\definecolor{darkblue}{rgb}{0.0, 0.0, 0.55}
\definecolor{bordeaux}{rgb}{0.34, 0.01, 0.1}
\usepackage[colorlinks,linkcolor=bordeaux,citecolor=darkblue,urlcolor=black,hypertexnames=true]{hyperref}
\usetikzlibrary{arrows,positioning}
\usepackage{colortbl}
\definecolor{myblue}{cmyk}{.3,0,0,0}
\definecolor{mypink}{cmyk}{0,.3,0,0}
\definecolor{mymix}{cmyk}{0.3,.3,0,0}
\usepackage{mathtools}

\newtheorem{theorem}{Theorem}[section]

\newtheorem{example}[theorem]{Example}

\numberwithin{equation}{section}

\def\R{{\mathbb{R}}}
\def\C{{\mathbb{C}}}
\def\N{{\mathbb{N}}}

\def\M{{\mathbf{M}}}
\def\x{{\mathbf{x}}}

\def\y{{\mathbf{y}}}

\def\o{{\mathbf{o}}}

\def\a{{\boldsymbol{\alpha}}}

\def\b{{\boldsymbol{\beta}}}
\def\g{{\boldsymbol{\gamma}}}

\def\A{{\mathscr{A}}}
\def\B{{\mathscr{B}}}
\def\CC{{\mathscr{C}}}

\def\S{{\mathbf{S}}}

\def\M{{\mathbf{M}}}

\def\supp{\hbox{\rm{supp}}}

\def\st{\hbox{\rm{s.t.}}}

\def\int{\hbox{\rm{int}}}

\newif\ifcomment
\commentfalse
\commenttrue
\usepackage{todonotes}

\newcommand{\revise}[1]{\textcolor{black}{#1}}
\begin{document}

\title[Certifying Global Optimality of AC-OPF Solutions]{Certifying Global Optimality of AC-OPF Solutions via sparse polynomial optimization}
% \author[J. Wang \and V. Magron \and J.B. Lasserre]{Jie Wang \and Victor Magron \and Jean B. Lasserre}

\author[J. Wang]{Jie Wang}
\address{Jie Wang: Academy of Mathematics and Systems Science, CAS}
\email{wangjie212@amss.ac.cn}
\urladdr{https://wangjie212.github.io/jiewang}

\author[V. Magron]{Victor Magron}
\address{Victor Magron: Laboratory for Analysis and Architecture of Systems, CNRS}
\email{vmagron@laas.fr}
\urladdr{https://homepages.laas.fr/vmagron}

\author[J. B. Lasserre]{Jean B. Lasserre}
\address{Jean B. Lasserre: Laboratory for Analysis and Architecture of Systems, CNRS}
\email{lasserre@laas.fr}
\urladdr{https://homepages.laas.fr/lasserre/}


\date{\today}

\begin{abstract}
We report the experimental results on certifying $1$\% global optimality of solutions of AC-OPF instances from PGLiB via the CS-TSSOS hierarchy --- a moment-SOS based hierarchy that exploits both correlative and term sparsity, which can provide tighter SDP relaxations than Shor's relaxation. Our numerical experiments demonstrate that the CS-TSSOS hierarchy scales well with the problem size and is indeed useful in certifying global optimality of solutions for large-scale real world problems, e.g., the AC-OPF problem. 
In particular, we are able to certify $1$\% global optimality for a challenging AC-OPF instance with $6515$ buses involving $14398$ real variables and $63577$ constraints.
\end{abstract}

\maketitle

\section{Introduction}

\paragraph{\textbf{Background on polynomial optimization}}
Polynomial optimization considers optimization problems where both the cost function and constraints are defined by polynomials, which widely arises in numerous fields, such as optimal power flow \cite{ghaddar2015optimal}, numerical analysis \cite{toms17}, computer vision \cite{yang2020one}, deep learning \cite{chen2020polynomial}, discrete optimization \cite{slot2020sum}, etc. Even though it is usually not difficult to find a locally optimal solution via a local solver (e.g., Ipopt \cite{wachter2006implementation}), the task of solving a polynomial optimization problem (POP) to global optimality is NP-hard in general. Over the last decades, the moment-sums of squares (moment-SOS) hierarchy consisting of a sequence of increasingly tight SDP relaxations, initially established by Lasserre \cite{Las01}, has become a popular tool to handle polynomial optimization. The moment-SOS hierarchy features its global convergence and finite convergence under mild conditions \cite{nie2014optimality}. However, the main concern on the moment-SOS hierarchy comes from its scalability as the $d$-th step of the moment-SOS hierarchy involves a semidefinite program (SDP) of size $\binom{n+d}{d}$ where $n$ is the number of decision variables of the POP. Except the first (relaxation) step of the moment-SOS hierarchy (also known as Shor's relaxation for quadratically constrained quadratic programs (QCQP) \cite{shor1987quadratic}), solving higher steps of the moment-SOS hierarchy is typically limited to small-scale POPs, at least when relying on interior-point solvers. To overcome this scalability issue, one practicable way is to exploit the structure of the POP to reduce the size of SDPs arising from the moment-SOS hierarchy. Such structures include symmetry \cite{riener2013exploiting}, correlative sparsity \cite{Las06,waki}, term sparsity \cite{tssos2,tssos1,tssos3}. The purpose of this paper is to demonstrate that the scalability of the moment-SOS hierarchy can be significantly improved when appropriate sparsity patterns are accessible via a thorough numerical experiment on the AC optimal power flow (AC-OPF) problem.

~

\paragraph{\textbf{Background on the AC-OPF problem}}
The AC-OPF is a fundamental problem in power systems, which has been extensively studied in recent years; for a detailed introduction and recent developments, the reader is referred to the survey \cite{bienstock2020} and references therein. One can formulate the AC-OPF problem as a POP either with real variables \cite{bienstock2020,ghaddar2015optimal} or with complex variables \cite{josz2018lasserre}. 
Nonlinear programming tools can mostly produce a locally optimal solution whose global optimality is however unknown. 
\revise{
We also mention that ill-conditioned power systems have proved the existence of computational pitfalls in which case the related nonlinear optimization algorithm can encounter unfeasibility issues,  especially in the radial parts of the networks. 
If we are looking precisely for such unfeasibility limit behaviors, one can often encounter this issue but rather rarely when one optimizes under realistic conditions, while taking into account the maximal injections of both producers and consumers.
In our study we are in the latter situation as we focus on real data coming from PGLiB instances. 
There have been several well-known heuristics plugged into local optimization schemes, e.g., Newton-Raphson, to make operational OPF problems more robust to solve \cite{tripathy1982load}. 
Alternatively, one can replace the Newton-Raphson scheme by a trust-region method (more accurate but also more CPU-demanding since information about the Hessian are required) as in \cite{abdelaziz2012globally} or a Riemannian optimization approach, which may yield better algorithmic convergence guarantees, as in \cite{heidarifar2021riemannian}. 
However such issues are less likely to arise when using semidefinite programming based on interior-point methods, since the underlying linear system to solve is different from the one considered in local solvers such as Newton-Raphson or trust-region schemes; a comprehensive study on this phenomenon could be found in \cite{eltved2019robustness}.}

Since 2006, several convex relaxation schemes (e.g., second order cone relaxations (SOCR) \cite{jabr2006radial}, quadratic convex relaxations (QCR) \cite{cof}, tight-and-cheap conic
relaxations (TCR) \cite{bingane2018tight} and semidefinite (Shor's) relaxations (SDR) \cite{bai2008semidefinite}) have been proposed to provide lower bounds for the AC-OPF which can be then used to certify global optimality of locally optimal solutions. \revise{In particular, SDR is equivalent to the first relaxation step of the moment-SOS hierarchy. SOCR is weaker but cheaper than SDR. It was shown that TCR offers a trade-off between SOCR and SDR in terms of optimality gap and computational cost \cite{bingane2018tight}, and that QCR is stronger than SOCR but neither dominates nor is dominated by SDR \cite{cof}. An experimental study in \cite{bingane2021conicopf} showed that on average TCR dominates QCR in terms of optimality gap.} While these relaxations (SOCR, QCR, TCR, SDR) could be scalable to problems of large size and prove to be tight for quite a few cases \cite{baba2019,cof,eltved2019robustness}, they yield significant optimality gaps for a large number of other cases\footnote{The reader may find related results on benchmarking SOCR and QCR at \href{https://github.com/power-grid-lib/pglib-opf/blob/master/BASELINE.md}{https://github.com/power-grid-lib/pglib-opf/blob/master/BASELINE.md}.}. \revise{To tackle these more challenging cases, it is then mandatory to go beyond SDR, providing tighter lower bounds. Along with this line recently in \cite{gopinath2020proving}, Gopinath et al. certified $1$\% global optimality for all AC-OPF instances with up to $300$ buses from the AC-OPF library \href{https://github.com/power-grid-lib/pglib-opf}{PGLiB} using an SDP-based bound tightening approach.
Relying on the complex moment-SOS hierarchy combined with a multi-order technique, Josz and Molzahn certified $0.05$\%
global optimality for certain $2000$-bus cases on a simplified AC-OPF model \cite{josz2018lasserre}.
To certify global optimality for challenging AC-OPF instances in the general form and with thousands of buses, the CS-TSSOS hierarchy then comes into play.}

~

\paragraph{\textbf{The CS-TSSOS hierarchy for large-scale POPs}}
The CS-TSSOS hierarchy developed in \cite{tssos3} by the authors is a sparsity-adapted version of the moment-SOS hierarchy targeted at large-scale POPs by simultaneously exploiting correlative sparsity (CS) and term sparsity (TS). The underlying idea is the following:
\begin{enumerate}
    \item[(1)] decomposing the system into subsystems by exploiting correlative sparsity, i.e., the fact that only a few variable products occur;
    \item[(2)] exploiting term sparsity, i.e., the fact that the input data only contain a few terms (by comparison with the maximal possible amount), to each subsystem to further reduce the size of SDPs.
\end{enumerate}
By virtue of this two-step reduction procedure, one may obtain SDP relaxations of significantly smaller size compared to the original SDP relaxations. \revise{Moreover, global convergence of the CS-TSSOS hierarchy can be still guaranteed under certain conditions \cite{tssos3}. Next the main concern on the CS-TSSOS hierarchy might be how it performs when being applied to real-word large-scale POPs in terms of scalability and accuracy, which will be addressed in the present paper.}

~

\paragraph{\textbf{Certifying global optimality for AC-OPF instances from PGLiB}}
\revise{We first propose a minimal initial relaxation step for the CS-TSSOS hierarchy. Being applied to the AC-OPF problem, this initial relaxation step is tighter than SDR and is less expensive than the second relaxation step of the CS-TSSOS hierarchy.} Then as the main contribution of this paper, we benchmark the CS-TSSOS hierarchy through a comprehensive numerical experiment on AC-OPF instances from the AC-OPF library \href{https://github.com/power-grid-lib/pglib-opf}{PGLiB} v20.07 \cite{baba2019} with up to tens of thousands of variables and constraints. The experimental results (see Section \ref{res}) demonstrate that the CS-TSSOS hierarchy scales well with the problem size and is able to certify global optimality (in the sense of reducing optimality gap within $1\%$) for many of the challenging test cases. In particular, the largest instance whose global optimality is certified beyond Shor's relaxation involves $14398$ real variables and $63577$ constraints. Besides, the largest instance for which the CS-TSSOS hierarchy is able to provide a smaller optimality gap than Shor's relaxation involves $24032$ real variables and $96805$ constraints. To the best of our knowledge, this is the first time in the literature that one can solve higher steps of the moment-SOS hierarchy other than Shor's relaxation for POPs of such large sizes.

~

\revise{
\begin{table}[htbp]
\caption{Notation.}\label{Notation}
\renewcommand\arraystretch{1.2}
\centering
\begin{tabular}{c|c}
\hline
\bf{Notation}&\bf{Meaning}\\
\hline
$\N$&$\{0,1,2,\ldots\}$\\
\hline
$\R$ (resp. $\C$) &the set of real (resp. complex) numbers\\
\hline
$[m]$&$\{1,2,\ldots,m\}$\\
\hline
$[m:n]$&$\{m,m+1,\ldots,n\}$\\
\hline
$\x=(x_1,\ldots,x_n)$&a tuple of variables\\
\hline
$\R[\x]$&the ring of real $n$-variate polynomials\\
\hline
$\N^n_d$&$\{(\alpha_i)_{i=1}^n\in\N^n\mid\sum_{i=1}^n\alpha_i\le d\}$\\
\hline
$\a,\b,\g$&vectors in $\N^n_d$\\
\hline
$|\cdot|$&the cardinality of a set\\
\hline
$\mathbf{S}^r$&the set of $r\times r$ symmetric matrices\\
\hline
$A\succeq0$&$A$ is a positive semidefinite (PSD) matrix\\
\hline
$\mathbf{S}_+^r$&the set of $r\times r$ PSD matrices\\
\hline
$A\circ B$&the Hadamard product of matrices $A$ and $B$\\
\hline
$G(V,E)$&a graph $G$ with nodes $V$ and edges $E$\\
\hline
$V(G)$ (resp. $E(G)$) &the node (resp. edge) set of a graph $G$\\
\hline
$G'$&a chordal extension of a graph $G$\\
\hline
$B_G$&the adjacency matrix of a graph $G$\\
\hline
\end{tabular}
\end{table}}

\section{Notation and preliminaries}
A polynomial $f\in\R[\x]$ can be written as $f(\x)=\sum_{\a\in\A}f_{\a}\x^{\a}$ with $\A\subseteq\N^n$, $f_{\a}\in\R$, and $\x^{\a}\coloneqq x_1^{\alpha_1}\cdots x_n^{\alpha_n}$. The \emph{support} of $f$ is defined by $\supp(f)\coloneqq \{\a\in\A\mid f_{\a}\ne0\}$. For $d\in\N$, the set $\x^{\N^n_{d}}\coloneqq (\x^{\a})_{\a\in\N^n_d}$ is called the {\em standard monomial basis} up to degree $d$. For convenience we abuse notation in the sequel, and denote by $\N^n_{d}$ instead of $\x^{\N^n_{d}}$ the standard monomial basis and use the exponent $\a$ to represent a monomial $\x^{\a}$. With $\y=(y_{\a})_{\a\in\N^n}\subseteq\R$ being a sequence indexed by $\N^n$, let $L_{\y}:\R[\x]\rightarrow\R$ be the linear functional $f=\sum_{\a}f_{\a}\x^{\a}\mapsto L_{\y}(f)=\sum_{\a}f_{\a}y_{\a}$.
For $\a\in\N^n,\A,\B\subseteq\N^n$,
let $\a+\B\coloneqq \{\a+\b\mid\b\in\B\}$ and $\A+\B\coloneqq \{\a+\b\mid\a\in\A,\b\in\B\}$. For $\b=(\beta_i)\in\N^n$, let $\supp(\b)\coloneqq \{i\in[n]\mid\beta_i\ne0\}$.

An {\em (undirected) graph} $G(V,E)$, or simply $G$, consists of a set of nodes $V$ and a set of edges $E\subseteq\{\{u,v\}\mid u\ne v,(u,v)\in V\times V\}$. The {\em adjacency matrix} of a graph $G$ is denoted by $B_G$ for which we put ones on its diagonal. A {\em clique} of a graph is a subset of nodes that induces a complete subgraph. A {\em maximal clique} is a clique that is not contained in any other clique. By definition, a {\em chordal graph} is a graph in which any cycle of length at least four has a chord\footnote{A chord is an edge that joins two nonconsecutive nodes in a cycle.}.
Any non-chordal graph $G(V,E)$ can be always extended to a chordal graph $G'(V,E')$ by adding appropriate edges to $E$, which is called a {\em chordal extension} of $G(V,E)$. The chordal extension of $G$ is usually not unique and the symbol $G'$ is used to represent any specific chordal extension of $G$ throughout the paper.

Given a graph $G(V,E)$, a symmetric matrix $Q$ with rows and columns indexed by $V$ is said to have sparsity pattern $G$ if $Q_{uv}=Q_{vu}=0$ whenever $u\ne v$ and $\{u,v\}\notin E$. Let $\mathbf{S}_G$ be the set of symmetric matrices with sparsity pattern $G$ and let $\Pi_{G}$ be the projection from $\mathbf{S}^{|V|}$ to the subspace $\mathbf{S}_G$, i.e., for $Q\in\mathbf{S}^{|V|}$,
\begin{equation}\label{sec2-eq7}
[\Pi_{G}(Q)]_{uv}=\begin{cases}
Q_{uv}, &\textrm{if }u=v\textrm{ or }\{u,v\}\in E,\\
0, &\textrm{otherwise}.
\end{cases}
\end{equation}
The set $\Pi_{G}(\mathbf{S}_+^{|V|})$ denotes matrices in $\mathbf{S}_G$ that have a PSD completion in the sense that
diagonal entries, and off-diagonal entries corresponding to edges of $G$ are fixed; other off-diagonal entries are free. More precisely, 
$\Pi_{G}(\mathbf{S}_+^{|V|})=\{\Pi_{G}(Q)\mid Q\in\mathbf{S}_+^{|V|}\}$.
For a chordal graph $G$, the following theorem due to Grone et al.~gives a characterization of matrices in the PSD completable cone $\Pi_{G}(\mathbf{S}_+^{|V|})$, which plays a crucial role in sparse semidefinite programming.
\begin{theorem}[\cite{grone1984}, Theorem 7]\label{sec2-thm2}
Let $G(V,E)$ be a chordal graph and assume that $C_1,\ldots,C_t$ are the list of maximal cliques of $G(V,E)$. Then a matrix $Q\in\Pi_{G}(\mathbf{S}_+^{|V|})$ if and only if $Q[C_i]\succeq0$ for $i=1,\ldots,t$, where $Q[C_i]$ denotes the principal submatrix of $Q$ indexed by the clique $C_i$.
\end{theorem}

\begin{figure}
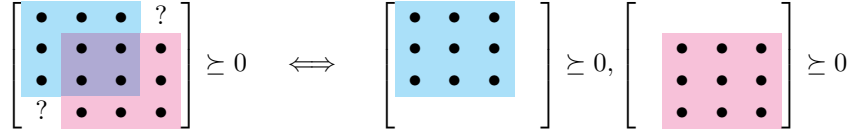

\centering
\begin{equation*}
    \begin{bmatrix}
    \begin{array}{cccc}
    \multicolumn{1}{>{\columncolor{myblue}}l}{\bullet}&\multicolumn{1}{>{\columncolor{myblue}}l}{\bullet}&\multicolumn{1}{>{\columncolor{myblue}}l}{\bullet}&?\\
    \multicolumn{1}{>{\columncolor{myblue}}l}{\bullet}&\multicolumn{1}{>{\columncolor{mymix}}l}{\bullet}&\multicolumn{1}{>{\columncolor{mymix}}l}{\bullet}&\multicolumn{1}{>{\columncolor{mypink}}l}{\bullet}\\
    \multicolumn{1}{>{\columncolor{myblue}}l}{\bullet}&\multicolumn{1}{>{\columncolor{mymix}}l}{\bullet}&\multicolumn{1}{>{\columncolor{mymix}}l}{\bullet}&\multicolumn{1}{>{\columncolor{mypink}}l}{\bullet}\\
    ?&\multicolumn{1}{>{\columncolor{mypink}}l}{\bullet}&\multicolumn{1}{>{\columncolor{mypink}}l}{\bullet}&\multicolumn{1}{>{\columncolor{mypink}}l}{\bullet}\\
    \end{array}
    \end{bmatrix}\succeq0
    \quad\iff\quad
    \begin{bmatrix}
    \begin{array}{ccccc}
    \multicolumn{1}{>{\columncolor{myblue}}l}{\bullet}&\multicolumn{1}{>{\columncolor{myblue}}l}{\bullet}&\multicolumn{1}{>{\columncolor{myblue}}l}{\bullet}&\\
    \multicolumn{1}{>{\columncolor{myblue}}l}{\bullet}&\multicolumn{1}{>{\columncolor{myblue}}l}{\bullet}&\multicolumn{1}{>{\columncolor{myblue}}l}{\bullet}&\\
    \multicolumn{1}{>{\columncolor{myblue}}l}{\bullet}&\multicolumn{1}{>{\columncolor{myblue}}l}{\bullet}&\multicolumn{1}{>{\columncolor{myblue}}l}{\bullet}&\\
    &&&\\
    \end{array}
    \end{bmatrix}\succeq0,
    \begin{bmatrix}
    \begin{array}{ccccc}
    &&&\\
    &\multicolumn{1}{>{\columncolor{mypink}}l}{\bullet}&\multicolumn{1}{>{\columncolor{mypink}}l}{\bullet}&\multicolumn{1}{>{\columncolor{mypink}}l}{\bullet}\\
    &\multicolumn{1}{>{\columncolor{mypink}}l}{\bullet}&\multicolumn{1}{>{\columncolor{mypink}}l}{\bullet}&\multicolumn{1}{>{\columncolor{mypink}}l}{\bullet}\\
    &\multicolumn{1}{>{\columncolor{mypink}}l}{\bullet}&\multicolumn{1}{>{\columncolor{mypink}}l}{\bullet}&\multicolumn{1}{>{\columncolor{mypink}}l}{\bullet}\\
    \end{array}
    \end{bmatrix}\succeq0
\end{equation*}
\caption{Illustration for Theorem \ref{sec2-thm2}.}
\end{figure}

\section{The CS-TSSOS hierarchy}
The moment-SOS hierarchy \cite{Las01} provides a sequence of increasingly tighter SDP relaxations for the following polynomial optimization problem:
\begin{equation}\label{pop}
(\textrm{POP}):\quad\begin{cases}
\inf\limits_{\x\in\R^n}&f(\x)\\ 
\,\,\,\textrm{s.t.}&g_j(\x)\ge 0,\quad j \in [m],\\
&g_{j}(\x)=0,\quad j\in[m+1:m+l],
\end{cases}
\end{equation}
where $f,g_1,\dots,g_{m+l}\in\R[\x]$ are all polynomials.

To state the moment hierarchy\footnote{We mainly focus on the moment hierarchy. The SOS hierarchy consists of the dual SDPs.}, recall that for a given $d\in\N$, the $d$-th order {\em moment matrix} $\M_{d}(\y)$ associated with $\y=(y_{\a})_{\a\in\N^n}$ is defined by $[\M_d(\y)]_{\b\g}\coloneqq L_{\y}(\x^{\b}\x^{\g})=y_{\b+\g},\forall\b,\g\in\N^n_{d}$ and the $d$-th order {\em localizing matrix} $\M_{d}(g\y)$ associated with $\y$ and $g=\sum_{\a}g_{\a}\x^{\a}\in\R[\x]$ is defined by
$[\M_{d}(g\,\y)]_{\b\g}\coloneqq L_{\y}(g\x^{\b}\x^{\g})=\sum_{\a}g_{\a}y_{\a+\b+\g},\forall\b,\g\in\N^n_{d}$. Let $d_j\coloneqq \lceil\deg(g_j)/2\rceil$ for $j=1,\ldots,m+l$ and $d_{\min}\coloneqq \max\,\{\lceil\deg(f)/2\rceil,d_1,\ldots,d_{m+l}\}$. Then for an integer $d\ge d_{\min}$, the $d$-th order moment relaxation for POP \eqref{pop} is given by
\begin{equation}\label{mom}
\begin{cases}
\inf\limits_{\y}&L_{\y}(f)\\
\textrm{s.t.}&\M_{d}(\y)\succeq0,\\
&\M_{d-d_j}(g_j\y)\succeq0,\quad j\in[m],\\
&\M_{d-d_j}(g_j\y)=0,\quad j\in[m+1:m+l],\\
&y_{\mathbf{0}}=1.
\end{cases}
\end{equation}

We call \eqref{mom} the {\em dense} moment hierarchy for POP \eqref{pop}, whose optima converge to the global optimum of \eqref{pop} under mild conditions (slightly stronger than compactness of the feasible set) \cite{Las01}. 
Unfortunately, when the relaxation order $d$ is greater than $1$, the dense moment hierarchy encounters a severe scalability issue as the maximal size of PSD constraints is a combinatorial number in terms of $n$ and $d$. Therefore in the following subsections, we briefly revisit the framework of exploiting sparsity to derive a {\em sparse} moment hierarchy of remarkably smaller size for POP \eqref{pop} in the presence of appropriate sparsity patterns. For details, the interested reader may refer to the early work on correlative sparsity by Waki et al. \cite{Las06,waki} and the recent work on term sparsity by the authors  \cite{wang2021exploiting,tssos2,tssos1,tssos3}.

\subsection{Correlative sparsity (CS)}\label{cs}
Let us from now on fix a relaxation order $d$. By exploiting correlative sparsity, we decompose the set of variables into a tuple of subsets and then the initial system splits into a tuple of subsystems. To this end, we define the {\em correlative sparsity pattern (csp) graph}\footnote{We adopt the idea of ``monomial sparsity" introduced in \cite{josz2018lasserre} for the definition of csp graphs, which thus is slightly different from the original definition given in \cite{waki}.} associated with POP \eqref{pop} to be the graph $G^{\textrm{csp}}$ with nodes $V=[n]$ and edges $E$ satisfying $\{i,j\}\in E$ if one of the following conditions holds:
\begin{enumerate}
    \item[(i)] there exists $\a\in\supp(f)\cup\bigcup_{k\in J'\cup K'}\supp(g_k)$ such that $\{i,j\}\subseteq\supp(\a)$;
    \item[(ii)] there exists $k\in[m+l]\setminus(J'\cup K')$ such that $\{i,j\}\subseteq\bigcup_{\a\in\supp(g_k)}\supp(\a)$,
\end{enumerate}
where $J'\coloneqq \{k\in[m]\mid d_k=d\}$ and $K'\coloneqq \{k\in[m+1:m+l]\mid d_k=d\}$.

Let $(G^{\textrm{csp}})'$ be a chordal extension of $G^{\textrm{csp}}$ and $\{I_k\}_{k\in[p]}$ be the list of maximal cliques of $(G^{\textrm{csp}})'$ with $n_k\coloneqq |I_k|$. We then partition the polynomials $g_j, j\in[m]\setminus J'$ into groups $\{g_j\mid j\in J_k\}, k\in[p]$ which satisfy
\begin{enumerate}
    \item[(i)] $J_1,\ldots,J_p\subseteq[m]\setminus J'$ are pairwise disjoint and $\cup_{k=1}^pJ_k=[m]\setminus J'$;
    \item[(ii)] for any $j\in J_k$, $\bigcup_{\a\in\supp(g_j)}\supp(\a)\subseteq I_k$, $k\in[p]$.
\end{enumerate}
Similarly, we also partition the polynomials $g_j, j\in[m+1:m+l]\setminus K'$ into groups $\{g_j\mid j\in K_k\}, k\in[p]$.

For any $k\in[p]$, let $\M_d(\y, I_k)$ (resp. $\M_d(g\y, I_k)$)
be the moment (resp. localizing) submatrix obtained from $\M_d(\y)$ (resp. $\M_d(g\y)$) by retaining only those rows and columns indexed by $\b\in\N_d^n$ of $\M_d(\y)$ (resp. $\M_d(g\y)$) with $\supp(\b)\subseteq I_k$.

\begin{example}\label{ex1}
Consider the POP:
\begin{equation*}
\begin{cases}
\inf\limits_{\x\in\R^3}&x_1^2x_2+x_2x_3^2\\ 
\,\,\st&1-x_1^2-x_2^2\ge0,\\
&1-x_2^2-x_3^2\ge0,\\
&x_1^4+x_2x_3=1.\\
\end{cases}
\end{equation*}
Let us take the relaxation order $d=2$. Then the csp graph is shown in Figure \ref{fg3}, which contains two maximal cliques: $\{x_1,x_2\}$ and $\{x_2,x_3\}$.

\begin{figure}[htbp]
\centering
\begin{tikzpicture}[every node/.style={circle, draw=blue!50, thick, minimum size=6mm}]
\node (n1) at (-1,-1.5) {$x_1$};
\node (n2) at (0,0) {$x_2$};
\node (n3) at (1,-1.5) {$x_3$};
\draw (n1)--(n2);
\draw (n2)--(n3);
\end{tikzpicture}
\caption{Illustration for correlative sparsity.}\label{fg3}
\end{figure}
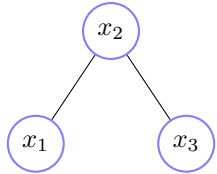
\end{example}

\subsection{Term sparsity (TS)}\label{cs-ts}
We next apply an iterative procedure to exploit term sparsity for each subsystem involving variables $\x(I_k)\coloneqq \{x_i\mid i\in I_k\}$ for $k\in[p]$. The intuition behind this procedure is the following: starting with a minimal initial set of moments, we expand the set of moments that is taken into account in the moment relaxation by iteratively performing chordal extension to the related graphs inspired by Theorem \ref{sec2-thm2}. More concretely, let $\mathscr{A}\coloneqq \supp(f)\cup\bigcup_{j=1}^{m+l}\supp(g_j)$ and $\mathscr{A}_k\coloneqq \{\a\in\A\mid\supp(\a)\subseteq I_k\}$ for $k\in[p]$. We define $G_{d,k,0}^{(0)}$ to be the graph with nodes $V_{d,k,0}\coloneqq \N^{n_k}_{d}$ and edges
\begin{equation}\label{ts-eq0}
E(G_{d,k,0}^{(0)})\coloneqq \left\{\{\b,\g\}\mid\b,\g\in V_{d,k,0},\b+\g\in\A_k\cup(2\N)^n\right\}.
\end{equation}
Note that here we embed $\N^{n_k}$ into $\N^{n}$ by specifying the $i$-th coordinate to be zero when $i\in[n]\setminus I_k$. 

\begin{example}
Consider again the POP in Example \ref{ex1} with the relaxation order $d=2$. Since there are two variable cliques derived from the correlative sparsity pattern, we have $p=2$. Figure \ref{fg4} illustrates the term sparsity pattern of this POP.

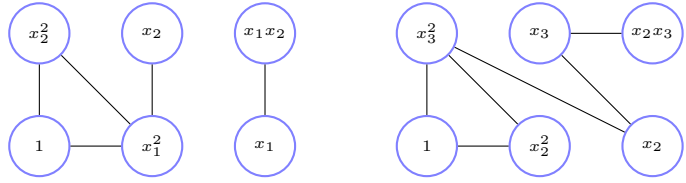
\begin{figure}[htbp]
{\tiny
\begin{minipage}{0.4\linewidth}
\centering
\begin{tikzpicture}[every node/.style={circle, draw=blue!50, thick, minimum size=8mm}]
\node (n1) at (0,0) {$1$};
\node (n2) at (1.5,0) {$x_1^2$};
\node (n3) at (3,0) {$x_1$};
\node (n4) at (0,1.5) {$x_2^2$};
\node (n5) at (1.5,1.5) {$x_2$};
\node (n6) at (3,1.5) {$x_1x_2$};
\draw (n1)--(n2);
\draw (n1)--(n4);
\draw (n2)--(n4);
\draw (n2)--(n5);
\draw (n3)--(n6);
\end{tikzpicture}
\end{minipage}
\begin{minipage}{0.4\linewidth}
\centering
\begin{tikzpicture}[every node/.style={circle, draw=blue!50, thick, minimum size=8mm}]
\node (n1) at (0,0) {$1$};
\node (n2) at (1.5,0) {$x_2^2$};
\node (n3) at (3,0) {$x_2$};
\node (n4) at (0,1.5) {$x_3^2$};
\node (n5) at (1.5,1.5) {$x_3$};
\node (n6) at (3,1.5) {$x_2x_3$};
\draw (n1)--(n2);
\draw (n1)--(n4);
\draw (n2)--(n4);
\draw (n3)--(n5);
\draw (n5)--(n6);
\draw (n3)--(n4);
\end{tikzpicture}
\end{minipage}}
\caption{Illustration for term sparsity: $G_{d,1,0}^{(0)}$ (left) and $G_{d,2,0}^{(0)}$ (right).}\label{fg4}
\end{figure}
\end{example}

For the sake of convenience, we set $g_0\coloneqq 1$ and $d_0\coloneqq 0$ hereafter and for a graph $G(V,E)$ with $V\subseteq\N^n$, let $\supp(G)\coloneqq\{\b+\g\mid\b=\g\in V\text{ or }\{\b,\g\}\in E\}$.
Assume that $G_{d,k,j}^{(0)},j\in J_k\cup K_k, k\in[p]$ are all empty graphs with $V_{d,k,j}\coloneqq \N^{n_k}_{d-d_j}$. Now for each $j\in\{0\}\cup J_k\cup K_k,k\in[p]$, we iteratively define an ascending chain of graphs $(G_{d,k,j}^{(s)})_{s\ge1}$ by
\begin{equation}\label{cts-eq3}
G_{d,k,j}^{(s)}\coloneqq (F_{d,k,j}^{(s)})',
\end{equation}
where $F_{d,k,j}^{(s)}$ is the graph with nodes $V_{d,k,j}$ and edges
\begin{equation}\label{cts-eq4}
E(F_{d,k,j}^{(s)})=\left\{\{\b,\g\}\mid\b,\g\in V_{d,k,j},(\b+\g+\supp(g_j))\cap\CC_{d}^{(s-1)}\ne\emptyset\right\},
\end{equation}
with
\begin{equation}\label{cts-eq2}
    \CC_{d}^{(s-1)}\coloneqq \bigcup_{k=1}^p\bigcup_{j\in \{0\}\cup J_k\cup K_k}\left(\supp(g_j)+\supp(G_{d,k,j}^{(s-1)})\right).
\end{equation}

Let $r_{d,k,j}\coloneqq |\N^{n_k}_{d-d_j}|=\binom{n_k+d-d_j}{d-d_j}$ for all $k,j$. Then for each $s\ge1$, the moment relaxation based on correlative-term sparsity for POP \eqref{pop} is given by
\begin{equation}\label{cts-eq5}
\begin{cases}
\inf\limits_{\y}&L_{\y}(f)\\
\textrm{s.t.}&B_{G_{d,k,0}^{(s)}}\circ \M_d(\y, I_k)\in\Pi_{G_{d,k,0}^{(s)}}(\S_+^{r_{d,k,0}}),\quad k\in[p],\\
&B_{G_{d,k,j}^{(s)}}\circ \M_{d-d_j}(g_j\y, I_k)\in\Pi_{G_{d,k,j}^{(k)}}(\S_+^{r_{d,k,j}}),\quad j\in J_k,k\in[p],\\
&B_{G_{d,k,j}^{(s)}}\circ \M_{d-d_j}(g_j\y, I_k)=0,\quad j\in K_k,k\in[p],\\
&L_{\y}(g_j)\ge0,\quad j\in J',\\
&L_{\y}(g_j)=0,\quad j\in K',\\
&y_{\mathbf{0}}=1.
\end{cases}
\end{equation}
The above hierarchy is called the {\em CS-TSSOS} hierarchy, which is indexed by two parameters: the relaxation order $d$ and the {\em sparse order} $s$.

\subsection{The minimal initial relaxation step}\label{initial}
For POP \eqref{pop}, suppose that $f$ is not a homogeneous polynomial or the polynomials $g_j,j\in[m+l]$ are of different degrees as in the case of the AC-OPF problem. Then instead of using the uniform minimum relaxation order $d_{\min}$, it might be more beneficial, from the computational point of view, to assign different relaxation orders to different subsystems obtained from the correlative sparsity pattern for the initial relaxation step of the CS-TSSOS hierarchy. To this end, we redefine the csp graph $G^{\textrm{icsp}}(V,E)$ as follows: $V=[n]$ and $\{i,j\}\in E$ whenever there exists $\a\in\supp(f)\cup\bigcup_{j\in[m+l]}\supp(g_j)$ such that $\{i,j\}\subseteq\supp(\a)$. This is clearly a subgraph of $G^{\textrm{csp}}$ defined in Section \ref{cs} and hence typically admits a smaller chordal extension. Let $(G^{\textrm{icsp}})'$ be a chordal extension of $G^{\textrm{icsp}}$ and $\{I_k\}_{k\in[p]}$ be the list of maximal cliques of $(G^{\textrm{icsp}})'$ with
$n_k\coloneqq |I_k|$. Now we partition the polynomials $g_j,j\in[m]$ into groups $\{g_j\mid j\in J_k\}_{k\in[p]}$ and $\{g_j\mid j\in J'\}$ which satisfy
\begin{enumerate}
    \item[(i)] $J_1,\ldots,J_p,J'\subseteq[m]$ are pairwise disjoint and $\bigcup_{k=1}^pJ_k\cup J'=[m]$;
    \item[(ii)] for any $j\in J_k$, $\bigcup_{\a\in\supp(g_j)}\supp(\a)\subseteq I_k$, $k\in[p]$;
     \item[(iii)] for any $j\in J'$, $\bigcup_{\a\in\supp(g_j)}\supp(\a)\nsubseteq I_k$ for all $k\in[p]$.
\end{enumerate}
Similarly, we also partition the polynomials $g_j,j\in[m+1:m+l]$ into groups $\{g_j\mid j\in K_k\}_{k\in[p]}$ and $\{g_j\mid j\in K'\}$.

Assume that $f$ decomposes as $f=\sum_{k\in[p]}f_k$ such that $\bigcup_{\a\in\supp(f_k)}\supp(\a)\subseteq I_k$ for $k\in[p]$. We define the vector of minimum relaxation orders $\o=(o_k)_k\in\N^{p}$ with $o_k\coloneqq \max\,(\{d_j:j\in J_k\cup K_k\}\cup\{\lceil\deg(f_k)/2\rceil\})$. Then with $s\ge1$, we define the following minimal initial relaxation step for the CS-TSSOS hierarchy:
\begin{equation}\label{cts}
\begin{cases}
\inf\limits_{\y} &L_{\y}(f)\\
\textrm{s.t.}&B_{G_{o_k,k,0}^{(s)}}\circ \M_{o_k}(\y, I_k)\in\Pi_{G_{o_k,k,0}^{(s)}}(\S_+^{t_{k,0}}),\quad k\in[p],\\
&\M_{1}(\y, I_k)\succeq0,\quad k\in[p],\\
&B_{G_{o_k,k,j}^{(s)}}\circ \M_{o_k-d_j}(g_j\y, I_k)\in\Pi_{G_{o_k,k,j}^{(s)}}(\S_+^{t_{k,j}}),\quad j\in J_k, k\in[p],\\
&L_{\y}(g_j)\ge0,\quad j\in J',\\
&B_{G_{o_k,k,j}^{(s)}}\circ \M_{o_k-d_j}(g_j\y, I_k)=0,\quad j\in K_k, k\in[p],\\
&L_{\y}(g_j)=0,\quad j\in K',\\
&y_{\mathbf{0}}=1,
\end{cases}
\end{equation}
where $G_{o_k,k,j}^{(s)},j\in J_k\cup K_k,k\in[p]$ are defined as in Section \ref{cs-ts} and $t_{k,j}\coloneqq \binom{n_k+o_k-d_j}{o_k-d_j}$ for all $k,j$. Note that in \eqref{cts} we add the PSD constraint on each first-order moment matrix $\M_{1}(\y, I_k)$ to strengthen the relaxation.
\vspace{1em}

The CS-TSSOS hierarchy is implemented in the Julia package {\tt TSSOS}\footnote{{\tt TSSOS} is freely available at \url{https://github.com/wangjie212/TSSOS}.}. In {\tt TSSOS}, the minimal initial relaxation step is accessible via the commands {\tt cs\_tssos\_first} and {\tt cs\_tssos\_higher!} by setting the relaxation order to be {\tt "min"}. For an introduction to {\tt TSSOS}, the reader is referred to \cite{magron2021tssos}.

\section{Problem formulation of AC-OPF}
The AC-OPF problem aims to minimize the generation cost of an alternating current transmission network under the physical constraints (Kirchhoff’s laws, Ohm’s law) as well as operational constraints, which can be formulated as the following POP in complex variables:
\begin{equation}\label{opf}
\begin{cases}
\inf\limits_{V_i,S_k^g\in\C}&\sum\limits_{k\in G}(\mathbf{c}_{2k}(\Re(S_{k}^g))^2+\mathbf{c}_{1k}\Re(S_{k}^g)+\mathbf{c}_{0k})\\
\quad\textrm{s.t.}&\angle V_r=0,\\
&\mathbf{S}_{k}^{gl}\le S_{k}^{g}\le \mathbf{S}_{k}^{gu},\quad\forall k\in G,\\
&\boldsymbol{\upsilon}_{i}^l\le|V_i|\le \boldsymbol{\upsilon}_{i}^u,\quad\forall i\in N,\\
&\sum_{k\in G_i}S_k^g-\mathbf{S}_i^d-\mathbf{Y}_i^s|V_{i}|^2=\sum_{(i,j)\in E_i\cup E_i^R}S_{ij},\quad\forall i\in N,\\
&S_{ij}=(\mathbf{Y}_{ij}^*-\mathbf{i}\frac{\mathbf{b}_{ij}^c}{2})\frac{|V_i|^2}{|\mathbf{T}_{ij}|^2}-\mathbf{Y}_{ij}^*\frac{V_iV_j^*}{\mathbf{T}_{ij}},\quad\forall (i,j)\in E,\\
&S_{ji}=(\mathbf{Y}_{ij}^*-\mathbf{i}\frac{\mathbf{b}_{ij}^c}{2})|V_j|^2-\mathbf{Y}_{ij}^*\frac{V_i^*V_j}{\mathbf{T}_{ij}^*},\quad\forall (i,j)\in E,\\
&|S_{ij}|\le \mathbf{s}_{ij}^u,\quad\forall (i,j)\in E\cup E^R,\\
&\boldsymbol{\theta}_{ij}^{\Delta l}\le \angle (V_i V_j^*)\le \boldsymbol{\theta}_{ij}^{\Delta u},\quad\forall (i,j)\in E.\\
\end{cases}
\end{equation}
The meaning of the symbols in \eqref{opf} is as follows: $N$ - the set of buses, $G$ - the set of generators, $G_i$ - the set of generators connected to bus $i$, $E$ - the set of {\em from} branches, $E^R$ - the set of {\em to} branches, $E_i$ and $E^R_i$ - the subsets of branches that are incident to bus $i$, $\mathbf{i}$ - imaginary unit, $V_i$ - the voltage at bus $i$, $S_k^g$ - the power generation at generator $k$, $S_{ij}$ - the power flow from bus $i$ to bus $j$, $\Re(\cdot)$ - real part of a complex number, $\angle(\cdot)$ - angle of a complex number, $|\cdot|$ - magnitude of a complex number, $(\cdot)^*$ - conjugate of a complex number, $r$ - the voltage angle reference bus. All symbols in boldface are constants ($\mathbf{c}_{0k},\mathbf{c}_{1k},\mathbf{c}_{2k},\boldsymbol{\upsilon}_{i}^l,\boldsymbol{\upsilon}_{i}^u,\mathbf{s}_{ij}^u,\boldsymbol{\theta}_{ij}^{\Delta l},\boldsymbol{\theta}_{ij}^{\Delta u}\in\R$, $\mathbf{S}_{k}^{gl},\mathbf{S}_{k}^{gu},\mathbf{S}_i^d,\mathbf{Y}_i^s,\mathbf{Y}_{ij},\mathbf{b}_{ij}^c,\mathbf{T}_{ij}\in\C$). For a full description on the AC-OPF problem, the reader may refer to \cite{baba2019}. By introducing real variables for both real and imaginary parts of each complex variable, we can convert the AC-OPF problem to a POP involving only real variables\footnote{The expressions involving angles of complex variables can be converted to polynomials by using $\tan(\angle z)=y/x$ for $z=x+\mathbf{i}y\in\C$.}.

To tackle an AC-OPF instance, we first compute a locally optimal solution with a local solver (e.g., {\tt Ipopt} \cite{wachter2006implementation}) and then rely on lower bounds obtained from certain relaxation schemes (SOCR/QCR/TCR/SDR/CS-TSSOS) to certify $1$\% global optimality. Suppose that the optimum reported by the local solver is AC and the lower bound given by a certain convex relaxation is opt. Then the {\em optimality gap} is defined by
\begin{equation*}
    \textrm{gap}\coloneqq \frac{\textrm{AC}-\textrm{opt}}{\textrm{AC}}\times100\%.
\end{equation*}
As in \cite{gopinath2020proving}, if the optimality gap is less than $1\%$, then we accept the locally optimal solution to be globally optimal. \revise{The procedure here for certifying global optimality of solutions of an AC-OPF instance is summarized in Algorithm \ref{sec5:alg3}.}

\begin{algorithm}
	\renewcommand{\algorithmicrequire}{\textbf{Input:}}
	\renewcommand{\algorithmicensure}{\textbf{Output:}}
	\caption{}\label{sec5:alg3}
	\begin{algorithmic}[1]
		\REQUIRE Data for an AC-OPF instance of \eqref{opf}
		\ENSURE A local solution and the optimality gap
		\STATE Build the POP model \eqref{pop} with the given data;
		\STATE Compute a local solution ($\coloneqq $ sol) of the POP model with a local solver, and denote the optimum by AC;
		\STATE Build the SDP relaxation \eqref{cts} for the POP model;
		\STATE Solve the SDP with {\tt Mosek}, and denote the optimum by opt;
		\STATE $\textrm{gap}\leftarrow\frac{\textrm{AC}-\textrm{opt}}{\textrm{AC}}\times100\%$;
		\RETURN sol, gap
	\end{algorithmic}
\end{algorithm}

\revise{In our experiments, \emph{we eliminate the power flow variables $S_{ij}$ from \eqref{opf} so that it only involves the voltage variable $V_i$ and the power generation variables $S_k^g$.} This reformulation is crucial to derive a CS-TSSOS hierarchy of lower complexity, which yet introduces a quartic constraint corresponding to $S_{ij}S^{*}_{ij}\le(\mathbf{s}_{ij}^u)^2$. \emph{To implement the first order relaxation, we then have to relax this quartic constraint to a quadratic constraint using the trick described in \cite[Sec. 5.3]{bienstock2020}.}} The minimal initial relaxation step of the CS-TSSOS hierarchy for \eqref{opf} is able to provide a tighter lower bound than the first order relaxation and is less expensive than the second order relaxation. Therefore, we hereafter refer to it as the 1.5th order relaxation. 

\section{Experimental settings}
\noindent{{\bf Challenging test cases.}}
Our test cases are selected from the AC-OPF library \href{https://github.com/power-grid-lib/pglib-opf}{PGLiB} v20.07 which provides various AC-OPF instances for benchmarking AC-OPF algorithms. For an introduction to this library, the reader is referred to \cite{baba2019}.
We observe that for a number of instances in PGLiB, the SOCR approach is able to close the gap to below $1\%$ and these instances are not particularly interesting as our purpose is to certify $1\%$ global optimality for more challenging cases. To that end, we select test cases from PGLiB (with no more than $25000$ buses) for which SOCR yields an optimality gap greater than $1\%$. There are $115$ such instances in total. For each instance, with {\tt TSSOS} we initially solve the first order relaxation and if this relaxation fails to certify $1$\% global optimality, we further solve the 1.5th order relaxation with $s=1$. Here {\tt Mosek 9.0} \cite{mosek} is employed as an SDP solver with the default settings.
\vspace{0.5em}

\noindent{{\bf Chordal extension.}}
To achieve a good balance between the computational cost and the approximation quality of lower bounds, two types of chordal extensions are used in the computation. For correlative sparsity, we use approximately smallest chordal extensions which give rise to small clique numbers. For term sparsity, instead we use maximal chordal extensions which make every connected component to be a complete subgraph by setting {\tt TS = "block"} in {\tt TSSOS}.
\vspace{0.5em}

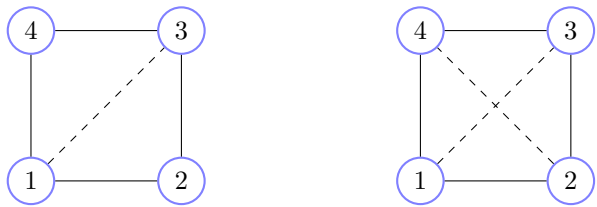
\begin{figure}[htbp]
\begin{minipage}{0.4\linewidth}
\centering
\begin{tikzpicture}[every node/.style={circle, draw=blue!50, thick, minimum size=6mm}]
\node (n1) at (0,0) {$1$};
\node (n2) at (2,0) {$2$};
\node (n4) at (0,2) {$4$};
\node (n3) at (2,2) {$3$};
\draw (n1)--(n2);
\draw (n2)--(n3);
\draw (n3)--(n4);
\draw (n4)--(n1);
\draw[dashed] (n1)--(n3);
\end{tikzpicture}
\end{minipage}
\begin{minipage}{0.4\linewidth}
\centering
\begin{tikzpicture}[every node/.style={circle, draw=blue!50, thick, minimum size=6mm}]
\node (n1) at (0,0) {$1$};
\node (n2) at (2,0) {$2$};
\node (n4) at (0,2) {$4$};
\node (n3) at (2,2) {$3$};
\draw (n1)--(n2);
\draw (n2)--(n3);
\draw (n3)--(n4);
\draw (n4)--(n1);
\draw[dashed] (n1)--(n3);
\draw[dashed] (n2)--(n4);
\end{tikzpicture}
\end{minipage}
\caption{Illustration for smallest chordal extension (left) and maximal chordal extension (right): the dashed edges are added via chordal extension.}\label{fg1}
\end{figure}

\noindent{{\bf Scaling of polynomial coefficients.}}
To improve the numerical conditioning of the SDP relaxations, we scale the coefficients of $f$ and $g_j$ so that they lie in the interval $[-1,1]$ before building the SDP relaxations.
\vspace{0.5em}

\noindent{{\bf Computational resources.}}
Instances with no more than $3500$ buses (except 2853\_sdet and 2869\_pegase) were computed on a laptop
with an Intel Core i5-8265U@1.60GHz CPU and 8GB RAM memory; instances with more than $3500$ buses (including 2853\_sdet and 2869\_pegase) were computed on a server with an Intel Xeon E5-2695v4@2.10GHz CPU and 128GB RAM memory.

\begin{table}[htbp]
\caption{Notations for the numerical results.}\label{not}
\renewcommand\arraystretch{1.2}
\centering
\begin{tabular}{c|c}
\hline
\bf{Notation}&\bf{Meaning}\\
\hline
AC&local optimum (available from PGLiB)\\
\hline
mc&maximal size of variable cliques\\
\hline
mb&maximal size of SDP blocks\\
\hline
opt&optimum of SDP relaxations\\
\hline
time&running time in seconds\\ 
\hline
gap&optimality gap ($\%$)\\
\hline
$*$&encountering a numerical error\\
\hline
-&out of memory\\
\hline
\end{tabular}
\end{table}

\section{Computational results and discussion}\label{res}
The computational results are summarized in Table \ref{opf-typ}--\ref{opf-sad} corresponding to three operational conditions, denoted by ``typical'', ``congested'' and ``small angle differences'', respectively, where the timing includes the time for pre-processing (to obtain the block structure), the time for building SDP and the time for solving SDP. Note that the maximal size of variable cliques varies from $6$ to $218$ among the tested cases. According to the tables, we can draw the following conclusions.
\vspace{0.5em}

\noindent{{\bf Reducing the optimality gap.}} As we would expect, the $1.5$th order relaxation provides tighter lower bounds than the first order relaxation. Indeed, when it is solvable, the $1.5$th order relaxation always reduces the optimality gap (unless the lower bounds given by the first order relaxation are already globally optimal). The largest instance for which the $1.5$th order relaxation is solvable is 10000\_goc and its corresponding POP involves $24032$ real variables and $96805$ constraints. The improvement of optimality gaps with the $1.5$th order relaxation is significant on quite a few cases. For instance, the first order relaxation yields a optimality gap of $42.96$\% on 30\_as under congested operating conditions while the $1.5$th order relaxation yields a optimality gap of merely $0.01$\%. 
On the other hand, as the cost of these improvements solving the $1.5$th order relaxation typically spends significantly more time than solving the first order relaxation.
\vspace{0.5em}

\noindent{{\bf Certifying $1$\% global optimality.}}
The first order relaxation is able to certify $1$\% global optimality for $29$ out of all $115$ instances. 
The $1.5$th order relaxation is able to certify $1$\% global optimality for $29$ out of the remaining $86$ instances. 
The largest instance for which we are able to certify $1$\% global optimality with the $1.5$th order relaxation is 6515\_rte and its corresponding POP involves $14398$ real variables and $63577$ constraints. One would expect that solving the second order relaxation could certify global optimality for more instances. 
However this is too expensive to implement for large-scale instances in practice.

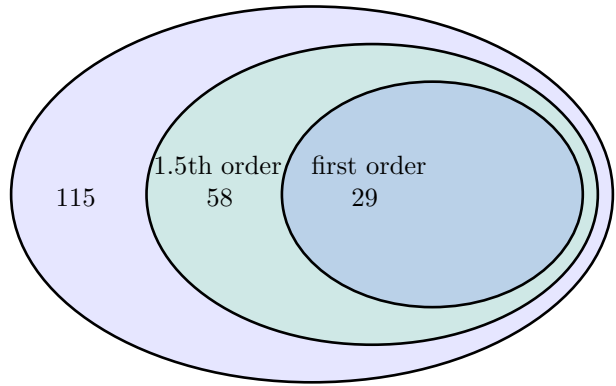
\begin{figure}[htbp]
\centering
\begin{tikzpicture}
\draw[draw=black,fill=blue,fill opacity=0.1,line width=1] (0,0) circle [x radius=4, y radius=2.5];
\node at (-2,0.2) {\parbox[c][2em][b]{8em}{115}};
\draw[draw=black,fill=green,fill opacity=0.1,line width=1] (0.8,0) circle [x radius=3, y radius=2];
\node at (-0.7,0.2) {\parbox[c][2em][b]{8em}{1.5th order\\ \hphantom{aaaa}58}};
\draw[draw=black,fill=blue,fill opacity=0.1,line width=1] (1.6,0) circle [x radius=2, y radius=1.5];
\node at (1.4,0.2) {\parbox[c][2em][b]{8em}{first order\\ \hphantom{aaa}29}};
\end{tikzpicture}
\caption{Certifying $1$\% global optimality for the test cases: the first order relaxation solves $29$ cases; the $1.5$th order relaxation solves extra $29$ cases.}\label{fg2}
\end{figure}

\vspace{0.5em}

\noindent{{\bf Computational burden.}}
The computational burden of the CS-TSSOS relaxations heavily relies on the maximal size of variable cliques. This is because large variable cliques usually lead to SDP matrices of large size in the resulting CS-TSSOS relaxations. It takes $63785$ seconds to solve the $1.5$th order relaxation for the case $4020$\_goc under congested operating conditions as it involves a variable clique of size 120. For similar reasons, {\tt Mosek} runs out of memory with the $1.5$th order relaxation for the cases $9241$\_pegase, $9591$\_goc, $10480$\_goc, $13659$\_pegase, $19402$\_goc, $24464$\_goc.
\vspace{0.5em}

\noindent{{\bf Numerical issues.}}
Even though we have scaled polynomial coefficients to improve numerical conditioning of the resulting SDPs, we observe that in numerous cases (especially when solving the $1.5$th order relaxation), the termination status of {\tt Mosek} is ``slow\_progress", which means that {\tt Mosek} does not converge to the default tolerance although the solver usually still gives a fairly good near-optimal solution in this case. 
Moreover, there are $12$ even more challenging instances for which {\tt Mosek} encounters severe numerical issues with the $1.5$th order relaxation and fails in converging to the optimum. 
This indicates that there is still room for improvement in order to tackle these challenging SDPs.

\begin{table}[htbp]
\caption{Results for AC-OPF problems: typical operating conditions.}\label{opf-typ}
\centering
\renewcommand\arraystretch{1.2}
\resizebox{\linewidth}{!}{
\begin{tabular}{|c|c|c|c|c|c|c|c|c|c|c|}
\hline
\multirow{2}*{case name}&\multirow{2}*{AC}&\multicolumn{4}{c|}{first order}&\multicolumn{5}{c|}{1.5th order}\\
\cline{3-11}
&&opt&time&mb&gap&mc&opt&time&mb&gap\\
\hline
3\_lmbd	&	5.8126e3	& 	5.7455e3	&0.10	&5	&1.15 	&6	&5.8126e3	&0.12	&22	&{\bf0.00}\\ 
\hline
5\_pjm	&	1.7552e4	& 	1.4997e4	&0.15	&6	&14.56 	&6	&1.7534e4	&0.58	&22	&{\bf0.10}\\
\hline
30\_ieee	&	8.2085e3	& 	7.5472e3	&0.22	&8	&8.06 	&8	&8.2085e3	&0.99	&22	&{\bf0.00} \\
\hline
% 39\_epri	&	1.3842e5	& 	1.3565e5	&0.68	&8	&2.00 	&8	&1.3842e5	&1.52	&25	&{\bf0.00} \\
% \hline
% 89\_pegase	&	1.0704e5	& 	1.0670e5	&1.18	&24	&{\bf0.32} 	&24	&1.0672e5	&1052	&184&{\bf0.30} \\
% \hline
% 118\_ieee	&	9.7214e4	& 	9.6900e4	&0.73	&10	&{\bf0.32} 	&10	&9.7214e4	&7.71	&37	&{\bf0.00} \\
% \hline
162\_ieee\_dtc&	1.0808e5	& 	1.0164e5	&2.15	&28	&5.96 	&28	&1.0645e5	&99.1	&74	&1.51 \\
\hline
240\_pserc	&	3.3297e6	& 	3.2512e6	&2.39	&16	&2.36 	&16	&3.3084e6	&28.6	&44	&{\bf0.64} \\
\hline
300\_ieee	&	5.6522e5	& 	5.5423e5	&2.72	&16	&1.94 	&14	&5.6522e5	&25.2	&40	&{\bf0.00} \\
\hline
588\_sdet	&	3.1314e5	& 	3.0886e5	&4.37	&18	&1.37 	&18	&3.1196e5	&50.6	&32	&{\bf0.38} \\
\hline
793\_goc	&	2.6020e5	& 	2.5636e5	&5.35	&18	&1.47 	&18	&2.5932e5	&66.1	&33	&{\bf0.34} \\
\hline
1888\_rte	&	1.4025e6	& 	1.3666e6	&30.0	    &26	&2.56 	&26	&1.3756e6	&458	&56	&1.92 \\
\hline
2312\_goc	&	4.4133e5	& 	4.3435e5	&87.8	&68	&1.58 	&68	&4.3858e5	&997	&81	&{\bf0.62} \\
\hline
2383wp\_k	&	1.8682e6	& 	1.8584e6	&63.0  	&50	&{\bf0.52} 	&48	&1.8646e6	&945	&77	&{\bf0.19} \\
\hline
2742\_goc	&	2.7571e5	& 	2.7561e5	&703	&92	&{\bf0.04} 	&   &			&	    &   &\\
\hline
% 2853\_sdet	&	2.0524e6	& 	2.0277e6	&114	&40	&1.20 	&40	&2.0493e6	&6711	&293&{\bf0.15}\\ 
% \hline
2869\_pegase&	2.4624e6	& 	2.4384e6	&85.0	    &26	&{\bf0.97} 	&26	&2.4571e6	&3641	&191&{\bf0.22} \\
\hline
3012wp\_k	&	2.6008e6	& 	2.5828e6	&123	&52	&{\bf0.69} 	&52	&2.5948e6	&1969	&81	&{\bf0.23} \\
\hline
3022\_goc	&	6.0138e5	& 	5.9277e5	&115	&48	&1.43 	&50	&5.9858e5	&1886	&76	&{\bf0.47} \\
\hline
% 3120sp\_k	&	2.1480e6	& 	2.1414e6	&139	&52	&{\bf0.31} 	&58	&2.1459e6	&1926	&70	&{\bf0.10} \\
% \hline
% 3375wp\_k	&	7.4382e6	& 	7.4103e6	&182	&58	&{\bf0.38} 	&54	&7.4312e6	&2741	&90	&{\bf0.09} \\
% \hline
4020\_goc	&	8.2225e5	& 	8.2208e5	&2356	&112&{\bf0.02}   &   & 			&		&   &\\
\hline
% 4601\_goc	&	8.2624e5	& 	8.2622e5	&1215	&108&{\bf0.00}   &	&			&	    &   &\\
% \hline
% 4619\_goc	&	4.7670e5	& 	4.7667e5	&1589	&82	&{\bf0.00} 	&	&			&       &   &\\
% \hline
4661\_sdet	&	2.2513e6	& 	2.2246e6	&25746	&204&1.18 	&218&$*$	        &$*$		&285&$*$\\	
\hline
4917\_goc	&	1.3878e6	& 	1.3658e6	&267	&64	&1.59 	&68	&1.3793e6	&29562	&110&{\bf0.61}\\ 
\hline
6468\_rte	&	2.0697e6	& 	2.0546e6	&415	&54	&{\bf0.73} 	&	&			&       &   &\\
\hline
6470\_rte	&	2.2376e6	& 	2.2060e6	&478	&54	&1.41 	&58	&$*$           &$*$      &98 &$*$\\
\hline
6495\_rte	&	3.0678e6	& 	2.6327e6	&426	&56	&14.18 	&54	&$*$           &$*$       &108&$*$\\
\hline
6515\_rte	&	2.8255e6	& 	2.6563e6	&460	&56	&5.99 	&54	&$*$           &$*$       &108&$*$\\
\hline
9241\_pegase&	6.2431e6	& 	6.1330e6	&982	&64	&1.76 &64 	&-			&-	    &1268   &-\\
\hline
% 9591\_goc	&	1.0617e6	& 		     	&   	&148&     &		&	        &       &   &\\
% \hline
10000\_goc	&	1.3540e6	& 	1.3460e6	&1714	&84	&{\bf0.59} &  	&			&       &   &\\
\hline
10480\_goc	&	2.3146e6	& 	2.3051e6	&8559	&136&{\bf0.41} &   &			&       &   &\\
\hline
13659\_pegase&	8.9480e6	& 	8.8707e6	&1808	&64	&{\bf0.86} &		&			&       &   &\\
\hline
19402\_goc	&	1.9778e6	& 	1.9752e6	&37157	&180&{\bf0.13} &		&			&       &   &\\
\hline
% 24464\_goc	&	2.6295e6	& 			    &	    &116&	  &	    &           &       &   &\\
% \hline
\end{tabular}}
\end{table}

\begin{table}[htbp]
\caption{Results for AC-OPF problems: congested operating conditions.}\label{opf-api}
\centering
\renewcommand\arraystretch{1.2}
\resizebox{\linewidth}{!}{
\begin{tabular}{|c|c|c|c|c|c|c|c|c|c|c|}
\hline
\multirow{2}*{case name}&\multirow{2}*{AC}&\multicolumn{4}{c|}{first order}&\multicolumn{5}{c|}{1.5th order}\\
\cline{3-11}
&&opt&time&mb&gap&mc&opt&time&mb&gap\\
\hline
3\_lmbd	    &1.1236e4	&1.0685e4	&0.11	&5	&4.90 	&6	&1.1236e4	&0.23	&22	&{\bf0.00} \\
\hline
5\_pjm	    &7.6377e4	&7.3253e4	&0.14	&6	&4.09 	&6	&7.6377e4	&0.55	&22	&{\bf0.00} \\
\hline
14\_ieee    &5.9994e3	&5.6886e3	&0.17	&6	&5.18 	&6	&5.9994e3	&0.54	&22	&{\bf0.00} \\
\hline
24\_ieee\_rts&1.3494e5	&1.2630e5	&0.37	&10	&6.40 	&10	&1.3392e5	&1.52	&31	&{\bf0.76} \\
\hline
30\_as	    &4.9962e3	&2.8499e3	&0.36	&8	&42.96 	&8	&4.9959e3	&2.41	&22	&{\bf0.01} \\
\hline
30\_ieee	&1.8044e4	&1.7253e4	&0.25	&8	&4.38 	&8	&1.8044e4	&1.24	&22	&{\bf0.00} \\
\hline
39\_epri	&2.4967e5	&2.4522e5	&0.28	&8	&1.78 	&8	&2.4966e5	&2.72	&25	&{\bf0.00} \\
\hline
73\_ieee\_rts&4.2263e5	&3.9912e5	&0.76	&12	&5.56 	&12	&4.1495e5	&6.77	&36	&1.82 \\
\hline
89\_pegase	&1.2781e5	&1.0052e5	&1.13	&24	&21.35 	&24	&1.0188e5	&1404	&184&20.29 \\
\hline
118\_ieee	&2.4224e5	&1.9375e5	&1.80	&10	&20.02 	&10	&2.2151e5	&11.5	&37	&8.56 \\
\hline
162\_ieee\_dtc&1.2099e5	&1.1206e5	&1.97	&28	&7.38 	&28	&1.1955e5	&84.1	&74	&1.19 \\
\hline
179\_goc	&1.9320e6	&1.7224e6	&1.24	&10	&10.85 	&10	&1.9226e6	&9.69	&37	&{\bf0.48} \\
\hline
% 240\_pserc	&4.6406e6	&4.6172e6	&2.28	&16	&{\bf0.50} 	&16	&4.6335e6	&27.9	&44	&{\bf0.15} \\
% \hline
% 300\_ieee	&6.8499e5	&6.7932e5	&2.08	&16	&{\bf0.83} 	&14	&6.8491e5	&31.6	&40	&{\bf0.01} \\
% \hline
500\_goc	&6.9241e5	&6.6004e5	&4.31	&18	&4.67 	&18	&6.7825e5	&78.0     &50	&2.05 \\
\hline
588\_sdet	&3.9476e5	&3.9026e5	&6.42	&18	&1.14 	&18	&3.9414e5	&57.0	    &32	&{\bf0.15} \\
\hline
793\_goc	&3.1885e5	&2.9796e5	&6.18	&18	&6.55 	&18	&3.1386e5	&79.2	&33	&1.56 \\
\hline
% 1354\_pegase&1.5135e6	&1.5025e6	&18.8	&26	&{\bf0.72} 	&26	&1.5122e6	&335	&49	&{\bf0.08} \\
% \hline
% 1951\_rte	&2.4108e6	&2.3635e6	&32.8	&26	&1.96 	&26	&2.4029e6	&596	&47	&{\bf0.32} \\
% \hline
2000\_goc	&1.4686e6	&1.4147e6	&54.0	&42	&3.67 	&42	&1.4610e6	&1094	&62	&{\bf0.51} \\
\hline
2312\_goc	&5.7152e5	&4.7872e5	&93.4	&68	&16.24 	&68	&5.2710e5	&972	&81	&7.77 \\
\hline
2736sp\_k	&6.5394e5	&5.8042e5	&89.6	&50	&11.24 	&48	&$*$           &$*$	    &79	&$*$\\
\hline
2737sop\_k	&3.6531e5	&3.4557e5	&71.9	&48	&5.40 	&48	&3.4557e5	&1653	&77	&5.40 \\
\hline
2742\_goc	&6.4219e5	&5.0824e5	&772	&92	&20.86 	&90	&6.0719e5	&4644	&108&5.45\\ 
\hline
2853\_sdet	&2.4578e6	&2.3869e6	&118	&40	&2.88 	&40	&2.4445e6	&10292	&293&{\bf0.54} \\
\hline
2869\_pegase&2.9858e6	&2.9604e6	&90.2	&26	&{\bf0.85} 	&26	&2.9753e6	&5409	&191&{\bf0.35} \\
\hline
3022\_goc	&6.5185e5	&6.2343e5	&102	&48	&4.36 	&50	&6.4070e5	&1519	&76	&1.71 \\
\hline
3120sp\_k	&9.3599e5	&7.6012e5	&138	&52	&18.79 	&58	&8.5245e5	&1627	&70	&8.93 \\
\hline
3375wp\_k	&5.8460e6	&5.5378e6	&222	&58	&5.27 	&54	&5.7148e6	&2619	&90	&2.25 \\
\hline
3970\_goc	&1.4241e6	&1.0087e6	&2469	&104&29.17 	&98	&1.0719e6	&15482	&135	&24.73 \\
\hline
4020\_goc	&1.2979e6	&1.0836e6	&3523	&112&16.51 	&120&1.1218e6	&63785	&174	&13.57\\ 
\hline
4601\_goc	&7.9253e5	&6.7523e5	&2143	&108&14.80 	&98	&7.3914e5	&17249	&125	&6.74 \\
\hline
4619\_goc	&1.0299e6	&9.6351e5	&1782	&82	&6.45   &84	&9.9766e5	&18348	&132	&3.13 \\
\hline
4661\_sdet	&2.6953e6	&2.6112e6	&15822	&204&3.12 	&218&$*$	        &$*$       &285	&$*$\\
\hline
4837\_goc	&1.1578e6	&1.0769e6	&500	&80	&6.98 	&84	&1.0947e6	&8723	&132	&5.45\\ 
\hline
4917\_goc	&1.5479e6	&1.4670e6	&259	&64	&5.23 	&68	&1.5180e6	&4688	&110	&1.93 \\
\hline
% 6468\_rte	&2.3135e6	&2.2920e6	&416	&54	&{\bf0.92} 	&	&			&       &       &\\
% \hline
6470\_rte	&2.6065e6	&2.5795e6	&427	&54	&1.04 	&58	&$*$           &$*$	    &	98	&$*$\\
\hline
6495\_rte	&2.9750e6	&2.9092e6	&453	&56	&2.21 	&54	&$*$           &$*$	    &	108	&$*$\\
\hline
6515\_rte   &3.0617e6	&2.9996e6	&421	&56	&2.02 	&54	&3.0434e6	&8456	&108	&{\bf0.60}\\ 
\hline
9241\_pegase&7.0112e6	&6.8784e6	&865	&64	&1.89 	&64	&-			&-      &1268   &-\\
\hline
9591\_goc	&1.4259e6	&1.2425e6   &7674	&148&12.86	&134&-			&-       &201    &-\\
\hline
10000\_goc	&2.3728e6	&2.1977e6	&2564	&84	&7.38 	&84	&2.3206e6	&27179	&97	    &2.20\\ 
\hline
10480\_goc	&2.7627e6	&2.6580e6	&8791   &136&3.79	&132&-			&-       &208    &-\\
\hline
13659\_pegase&9.2842e6	&9.1360e6	&1599	&64	&1.60 	&64	&-			&-       &1268   &-\\
\hline
19402\_goc	&2.3987e6	&2.3290e6	&32465  &180&2.91	&172&-	        &-       &242   &-\\
\hline
24464\_goc	&2.4723e6	&2.4177e6	&11760  &116&2.21	&118&-	        &-       &172   &-\\
\hline
\end{tabular}}
\end{table}

\begin{table}[htbp]
\caption{Results for AC-OPF problems: small angle difference conditions.}\label{opf-sad}
\centering
\renewcommand\arraystretch{1.2}
\resizebox{\linewidth}{!}{
\begin{tabular}{|c|c|c|c|c|c|c|c|c|c|c|}
\hline
\multirow{2}*{case name}&\multirow{2}*{AC}&\multicolumn{4}{c|}{first order}&\multicolumn{5}{c|}{1.5th order}\\
\cline{3-11}
&&opt&time&mb&gap&mc&opt&time&mb&gap\\
\hline
3\_lmbd	&5.9593e3	&5.7463e3	&0.11	&5	&3.57 	&6	&5.9593e3	&0.13	&22	&{\bf0.00} \\
\hline
5\_pjm	&2.6109e4	&2.6109e4	&0.12	&6	&{\bf0.00} 	&	&			&       &   &\\
\hline
14\_ieee&2.7768e3	&2.7743e3	&0.14	&6	&{\bf0.09} 	&	&			&       &   &\\
\hline
24\_ieee\_rts&7.6918e4&7.3555e4	&0.21	&10	&4.37 	&10	&7.4852e4	&1.64	&31	&2.69 \\
\hline
30\_as	&8.9735e2	&8.9527e2	&0.19	&8	&{\bf0.23} 	&	&			&       &   &\\
\hline
30\_ieee	&8.2085e3	&7.5472e3	&0.30	&8	&8.06 	&8	&8.2085e3	&1.03	&22	&{\bf0.00} \\
\hline
% 39\_epri	&1.4834e5   &1.4791e5	&0.31	&8	&{\bf0.29} 	&8	&1.4831e5	&1.45	&25	&{\bf0.02} \\
% \hline
% 57\_ieee	&3.8663e4	&3.8646e4	&0.44	&12	&{\bf0.04} 	&	&			&       &   &\\
% \hline
73\_ieee\_rts&2.2760e5	&2.2136e5	&0.55	&12	&2.74 	&12	&2.2447e5	&5.01	&36	&1.38 \\
\hline
% 89\_pegase	&1.0729e5	&1.0671e5	&1.15	&24	&{\bf0.54} 	&24	&1.0672e5	&1416	&184&{\bf0.54} \\
% \hline
118\_ieee	&1.0516e5	&1.0191e5	&0.79	&10	&3.10 	&10	&1.0313e5	&10.8	&37	&1.93 \\
\hline
162\_ieee\_dtc&1.0869e5	&1.0282e5	&2.46	&28	&5.40 	&28	&1.0740e5	&105	&74	&1.19 \\
\hline
179\_goc	&7.6186e5	&7.5261e5	&1.39	&10	&1.21 	&10	&7.5573e5	&11.8	&37	&{\bf0.80} \\
\hline
240\_pserc	&3.4054e6	&3.2772e6	&3.00	&16	&3.76 	&16	&3.3128e6	&33.8	&44	&2.72 \\
\hline
300\_ieee	&5.6570e5	&5.6162e5	&2.73	&16	&{\bf0.72} 	&14	&5.6570e5	&25.2	&40	&{\bf0.00} \\
\hline
500\_goc	&4.8740e5	&4.6043e5	&6.52	&18	&5.53 	&18	&4.6098e5	&67.7	&50	&5.42 \\
\hline
588\_sdet	&3.2936e5	&3.1233e5	&5.12	&18	&5.17 	&18	&3.1898e5	&56.6	&32	&3.15 \\
\hline
793\_goc	&2.8580e5	&2.7133e5	&5.61	&18	&5.06 	&18	&2.7727e5	&76.0	&33	&2.98 \\
\hline
1354\_pegase&1.2588e6	&1.2172e6	&19.8	&26	&3.31 	&26	&1.2582e6	&387	&49	&{\bf0.05} \\
\hline
1888\_rte	&1.4139e6	&1.3666e6	&31.2	&26	&3.34 	&26	&1.3756e6	&497	&56	&2.71 \\
\hline
2000\_goc	&9.9288e5	&9.8400e5	&50.9	&42	&{\bf0.89} 	&42	&9.8435e5	&1052	&62	&{\bf0.86} \\
\hline
2312\_goc	&4.6235e5	&4.4719e5	&121	&68	&3.28 	&68	&4.5676e5	&1009	&81	&1.21 \\
\hline
2383wp\_k	&1.9112e6	&1.9041e6	&65.6	&50	&{\bf0.37} 	&48	&1.9060e6	&937	&77	&{\bf0.27} \\
\hline
2736sp\_k	&1.3266e6	&1.3229e6	&89.5	&50	&{\bf0.28} 	&	&			&       &   &\\
\hline
2737sop\_k	&7.9095e5	&7.8672e5	&76.3	&48	&{\bf0.53} 	&	&			&       &   &\\
\hline
2742\_goc	&2.7571e5	&2.7561e5	&686	&92	&{\bf0.04} 	&	&			&       &   &\\
\hline
2746wop\_k	&1.2337e6	&1.2248e6	&79.1	&48	&{\bf0.72} 	&	&			&       &   &\\
\hline
2746wp\_k	&1.6669e6	&1.6601e6	&83.1	&50	&{\bf0.41} 	&	&			&       &   & \\
\hline
2853\_sdet	&2.0692e6	&2.0303e6	&106	&40	&1.88 	&40	&2.0537e6	&40671	&293&{\bf0.75} \\
\hline
% 2868\_rte	&2.0213e6	&2.0078e6	&70	    &36	&{\bf0.67} 	&	&			&       &   &\\
% \hline
2869\_pegase&2.4687e6	&2.4477e6	&85.4	&26	&{\bf0.85} 	&	&			&       &   &\\
\hline
3012wp\_k	&2.6192e6	&2.5994e6	&97.1	&52	&{\bf0.76} 	&	&			&       &   &\\
\hline
3022\_goc	&6.0143e5	&5.9278e5	&93.4	&48	&1.44 	&50	&5.9859e5	&1340	&76	&{\bf0.47} \\
\hline
3120sp\_k	&2.1749e6	&2.1611e6	&117	&52	&{\bf0.64 }	&	&			&       &   &\\
\hline
% 3375wp\_k	&7.4382e6	&7.4103e6	&151	&58	&{\bf0.38} 	&	&			&       &   &\\
% \hline
% 3970\_goc	&9.6555e5	&9.6393e5	&1651	&104&{\bf0.17} 	&	&			&       &   &\\
% \hline
4020\_goc	&8.8969e5	&8.4238e5	&2746	&112&5.32 	&120&8.7038e5	&43180	&174&2.17 \\
\hline
4601\_goc	&8.7803e5	&8.3370e5	&1763	&108&5.05 	&98	&8.3447e5	&15585	&125&4.96 \\
\hline
4619\_goc	&4.8435e5	&4.8106e5	&1387	&82	&{\bf0.68} 	&	&			&       &   &\\
\hline
4661\_sdet	&2.2610e6	&2.2337e6	&16144	&204&1.21 	&218&$*$	        &$*$     	&285&$*$\\
\hline
% 4837\_goc	&8.7712e5	&8.7447e5	&429	&80	&{\bf0.30} 	&	&			&       &   &\\
% \hline
4917\_goc	&1.3890e6	&1.3665e6	&260	&64	&1.62 	&68	&1.3800e6	&4914	&110&{\bf0.65} \\
\hline
6468\_rte	&2.0697e6	&2.0546e6	&399	&54	&{\bf0.73} 	&	&			&       &   &\\
\hline
6470\_rte	&2.2416e6	&2.2100e6	&451	&54	&1.41 	&58	&$*$           &$*$	    &98	&$*$\\
\hline
6495\_rte	&3.0678e6	&2.6323e6	&404	&56	&14.19 	&54	&$*$           &$*$	    &108&$*$\\
\hline
6515\_rte	&2.8698e6	&2.6565e6	&399	&56	&7.43 	&54	&$*$           &$*$	    &108&$*$\\
\hline
9241\_pegase&6.3185e6	&6.1696e6	&912	&64	&2.36 	&64	&-			&-      &1268&-\\
\hline
9591\_goc	&1.1674e6	&1.0712e6	&7835	&148&8.24 	&134&-			&-       &201&-\\
\hline
10000\_goc	&1.4902e6	&1.4204e6	&2508	&84	&4.68 	&84	&1.4212e6	&25340	&97	&4.63 \\
\hline
10480\_goc	&2.3147e6	&2.3051e6	&6522   &136&{\bf0.42}	&&		    &      && \\
\hline
13659\_pegase&9.0422e6	&8.9142e6	&1653	&64	&1.42 	&64	&-			&-      &1268&-\\
\hline
19402\_goc	&1.9838e6	&1.9783e6	&30122  &180&{\bf0.28}	&&		    &      & &\\
\hline
24464\_goc	&2.6540e6   &2.6268e6	&12101  &116&1.03		&118&-	    	&-       &172 &-\\
\hline
\end{tabular}}
\end{table}

\section{Conclusions}
We have benchmarked the CS-TSSOS hierarchy on a number of challenging AC-OPF cases and demonstrated that the $1.5$th order relaxation is indeed useful in reducing the optimality gap and certifying global optimality of AC-OPF solutions. We hope that the computational results presented in this paper would convince people to think of CS-TSSOS as an alternative tool for certifying global optimality of solutions of large-scale POPs. One line of future research is to improve the efficiency of the CS-TSSOS relaxations by relying on more advanced chordal extension algorithms.
We also plan to design suitable branch and bound algorithms to reach better accuracy results such as $0.1$\% or $0.01$\% global optimality for the AC-OPF problem.
\revise{
Eventually, such a current degree of optimality, i.e., $1\%$ cannot be judged without considering uncertainty. 
For power systems of any size, the uncertainties of loads, generations, storage as well as configurations, impact upon steady-state operating points.
In a further dedicated study, we shall extend our presented framework to handle uncertainties by means of sparsity-adapted versions of the existing approaches for robust/parametric polynomial optimization \cite{laraki2012semidefinite,lasserre2010joint}. }

~

\paragraph{\textbf{Acknowledgements}.} 
We would like to thank Jean Maeght, Patrick Panciatici and Manual Ruiz for interesting discussions and insights regarding AC-OPF.
Both authors were supported by the Tremplin ERC Stg Grant ANR-18-ERC2-0004-01 (T-COPS project).
The second author was supported by the FMJH Program PGMO (EPICS project), as well as the PEPS2 Program (FastOPF project) funded by AMIES and RTE.
This work has benefited from  the European Union's Horizon 2020 research and innovation program under the Marie Sklodowska-Curie Actions, grant agreement 813211 (POEMA) as well as from the AI Interdisciplinary Institute ANITI funding, through the French ``Investing for the Future PIA3'' program under the Grant agreement n$^{\circ}$ANR-19-PI3A-0004.

\bibliographystyle{plain}
\bibliography{refer}

\begin{thebibliography}{10}

\bibitem{abdelaziz2012globally}
Morad~MA Abdelaziz, Hany~E Farag, EF~El-Saadany, and YA-RI Mohamed.
\newblock A globally convergent trust-region method for power flow studies in
  active distribution systems.
\newblock In {\em 2012 IEEE Power and Energy Society General Meeting}, pages
  1--7. IEEE, 2012.

\bibitem{mosek}
Mosek ApS.
\newblock The {MOSEK} optimization software, 2019.

\bibitem{baba2019}
Sogol Babaeinejadsarookolaee, Adam Birchfield, Richard~D Christie, Carleton
  Coffrin, Christopher DeMarco, Ruisheng Diao, Michael Ferris, Stephane
  Fliscounakis, Scott Greene, Renke Huang, et~al.
\newblock The power grid library for benchmarking {AC} optimal power flow
  algorithms.
\newblock {\em arXiv preprint arXiv:1908.02788}, 2019.

\bibitem{bai2008semidefinite}
Xiaoqing Bai, Hua Wei, Katsuki Fujisawa, and Yong Wang.
\newblock Semidefinite programming for optimal power flow problems.
\newblock {\em International Journal of Electrical Power \& Energy Systems},
  30(6-7):383--392, 2008.

\bibitem{bienstock2020}
Dan Bienstock, Mauro Escobar, Claudio Gentile, and Leo Liberti.
\newblock Mathematical programming formulations for the alternating current
  optimal power flow problem.
\newblock {\em 4OR}, 18(3):249--292, 2020.

\bibitem{bingane2018tight}
Christian Bingane, Miguel~F Anjos, and S{\'e}bastien Le~Digabel.
\newblock Tight-and-cheap conic relaxation for the ac optimal power flow
  problem.
\newblock {\em IEEE Transactions on Power Systems}, 33(6):7181--7188, 2018.

\bibitem{bingane2021conicopf}
Christian Bingane, Miguel~F Anjos, and S{\'e}bastien Le~Digabel.
\newblock Conicopf: Conic relaxations for ac optimal power flow computations.
\newblock In {\em 2021 IEEE Power \& Energy Society General Meeting (PESGM)},
  pages 1--5. IEEE, 2021.

\bibitem{chen2020polynomial}
Tong Chen, Jean~B Lasserre, Victor Magron, and Edouard Pauwels.
\newblock {Semialgebraic Optimization for Lipschitz Constants of ReLU
  Networks}.
\newblock In H.~Larochelle, M.~Ranzato, R.~Hadsell, M.~F. Balcan, and H.~Lin,
  editors, {\em Advances in Neural Information Processing Systems}, volume~33,
  pages 19189--19200. Curran Associates, Inc., 2020.

\bibitem{cof}
Carleton Coffrin, Hassan~L Hijazi, and Pascal Van~Hentenryck.
\newblock The {QC} relaxation: A theoretical and computational study on optimal
  power flow.
\newblock {\em IEEE Transactions on Power Systems}, 31(4):3008--3018, 2015.

\bibitem{eltved2019robustness}
Anders Eltved, Joachim Dahl, and Martin~S Andersen.
\newblock On the robustness and scalability of semidefinite relaxation for
  optimal power flow problems.
\newblock {\em Optimization and Engineering}, pages 1--18, 2019.

\bibitem{ghaddar2015optimal}
Bissan Ghaddar, Jakub Marecek, and Martin Mevissen.
\newblock Optimal power flow as a polynomial optimization problem.
\newblock {\em IEEE Transactions on Power Systems}, 31(1):539--546, 2015.

\bibitem{gopinath2020proving}
S~Gopinath, Hassan~L Hijazi, Tillmann Weisser, Harsha Nagarajan, Mertcan
  Yetkin, Kaarthik Sundar, and Russel~W Bent.
\newblock Proving global optimality of {ACOPF} solutions.
\newblock {\em Electric Power Systems Research}, 189, 2020.

\bibitem{grone1984}
Robert Grone, Charles~R Johnson, Eduardo~M S{\'a}, and Henry Wolkowicz.
\newblock Positive definite completions of partial hermitian matrices.
\newblock {\em Linear algebra and its applications}, 58:109--124, 1984.

\bibitem{heidarifar2021riemannian}
Majid Heidarifar, Panagiotis Andrianesis, and Michael Caramanis.
\newblock A riemannian optimization approach to the radial distribution network
  load flow problem.
\newblock {\em Automatica}, 129:109620, 2021.

\bibitem{jabr2006radial}
Rabih~A Jabr.
\newblock Radial distribution load flow using conic programming.
\newblock {\em IEEE transactions on power systems}, 21(3):1458--1459, 2006.

\bibitem{josz2018lasserre}
C{\'e}dric Josz and Daniel~K Molzahn.
\newblock Lasserre hierarchy for large scale polynomial optimization in real
  and complex variables.
\newblock {\em SIAM Journal on Optimization}, 28(2):1017--1048, 2018.

\bibitem{laraki2012semidefinite}
Rida Laraki and Jean~B Lasserre.
\newblock Semidefinite programming for min--max problems and games.
\newblock {\em Mathematical programming}, 131(1):305--332, 2012.

\bibitem{Las01}
J.-B. Lasserre.
\newblock {Global Optimization with Polynomials and the Problem of Moments}.
\newblock {\em SIAM Journal on Optimization}, 11(3):796--817, 2001.

\bibitem{Las06}
J.-B. Lasserre.
\newblock Convergent {SDP}-relaxations in polynomial optimization with
  sparsity.
\newblock {\em SIAM Journal on Optimization}, 17(3):822--843, 2006.

\bibitem{lasserre2010joint}
Jean~B Lasserre.
\newblock A ``joint+ marginal'' approach to parametric polynomial optimization.
\newblock {\em SIAM Journal on Optimization}, 20(4):1995--2022, 2010.

\bibitem{toms17}
V.~Magron, G.~Constantinides, and A.~Donaldson.
\newblock {Certified Roundoff Error Bounds Using Semidefinite Programming}.
\newblock {\em ACM Trans. Math. Softw.}, 43(4):1--34, 2017.

\bibitem{magron2021tssos}
Victor Magron and Jie Wang.
\newblock {TSSOS}: a {J}ulia library to exploit sparsity for large-scale
  polynomial optimization.
\newblock {\em arXiv preprint arXiv:2103.00915}, 2021.

\bibitem{nie2014optimality}
J.~Nie.
\newblock Optimality conditions and finite convergence of {L}asserre's
  hierarchy.
\newblock {\em Mathematical Programming}, 146(1):97--121, Aug 2014.

\bibitem{riener2013exploiting}
Cordian Riener, Thorsten Theobald, Lina~Jansson Andr{\'e}n, and Jean~B
  Lasserre.
\newblock Exploiting symmetries in {SDP}-relaxations for polynomial
  optimization.
\newblock {\em Mathematics of Operations Research}, 38(1):122--141, 2013.

\bibitem{shor1987quadratic}
Naum~Z Shor.
\newblock Quadratic optimization problems.
\newblock {\em Soviet Journal of Computer and Systems Sciences}, 25:1--11,
  1987.

\bibitem{slot2020sum}
Lucas Slot and Monique Laurent.
\newblock Sum-of-squares hierarchies for binary polynomial optimization.
\newblock {\em arXiv preprint arXiv:2011.04027}, 2020.

\bibitem{tripathy1982load}
SC~Tripathy, G~Durga Prasad, OP~Malik, and GS~Hope.
\newblock Load-flow solutions for ill-conditioned power systems by a
  newton-like method.
\newblock {\em IEEE Transactions on Power apparatus and Systems},
  (10):3648--3657, 1982.

\bibitem{wachter2006implementation}
Andreas W{\"a}chter and Lorenz~T Biegler.
\newblock On the implementation of an interior-point filter line-search
  algorithm for large-scale nonlinear programming.
\newblock {\em Mathematical programming}, 106(1):25--57, 2006.

\bibitem{waki}
H.~Waki, S.~Kim, M.~Kojima, and M.~Muramatsu.
\newblock {Sums of Squares and Semidefinite Programming Relaxations for
  Polynomial Optimization Problems with Structured Sparsity}.
\newblock {\em SIAM Journal on Optimization}, 17(1):218--242, 2006.

\bibitem{wang2021exploiting}
Jie Wang and Victor Magron.
\newblock Exploiting sparsity in complex polynomial optimization.
\newblock {\em arXiv preprint arXiv:2103.12444}, 2021.

\bibitem{tssos2}
Jie Wang, Victor Magron, and Jean-Bernard Lasserre.
\newblock {Chordal-{TSSOS}: a moment-{SOS} hierarchy that exploits term
  sparsity with chordal extension}.
\newblock {\em SIAM Journal on Optimization}, 2020.
\newblock Accepted for publication.

\bibitem{tssos1}
Jie Wang, Victor Magron, and Jean-Bernard Lasserre.
\newblock {TSSOS}: A moment-{SOS} hierarchy that exploits term sparsity.
\newblock {\em SIAM Journal on Optimization}, 31(1):30--58, 2021.

\bibitem{tssos3}
Jie Wang, Victor Magron, Jean-Bernard Lasserre, and Ngoc Hoang~Anh Mai.
\newblock {{CS-TSSOS}: Correlative and term sparsity for large-scale polynomial
  optimization}.
\newblock {\em arXiv:2005.02828}, 2020.

\bibitem{yang2020one}
Heng Yang and Luca Carlone.
\newblock One ring to rule them all: Certifiably robust geometric perception
  with outliers.
\newblock {\em arXiv preprint arXiv:2006.06769}, 2020.

\end{thebibliography}
\end{document}